\input amstex
\documentstyle {amsppt}
\UseAMSsymbols \vsize 18cm

\magnification 1200

\catcode`\@=11
\def\displaylinesno #1{\displ@y\halign{
\hbox to\displaywidth{$\@lign\hfil\displaystyle##\hfil$}&
\llap{$##$}\crcr#1\crcr}}

\def\ldisplaylinesno #1{\displ@y\halign{
\hbox to\displaywidth{$\@lign\hfil\displaystyle##\hfil$}&
\kern-\displaywidth\rlap{$##$} \tabskip\displaywidth\crcr#1\crcr}}
\catcode`\@=12

\refstyle{C}

\let \ti=\widetilde
  \font
\nr=eufb7 at 10pt  \font \srm=cmr10 at
7.5pt  \font\ens=msbm10
 \font\main=cmsy10 at 10pt
    \font \fin=lasy8 at
15.4 pt
\def \X{\mathop{\hbox{\nr X}^{\hbox{\srm nr}}}\nolimits}

\def \Xim{\mathop{\hbox{\nr X}^{\hbox{\srm nr}}_0}\nolimits}

\def \o{\mathop{\hbox{\main O}}\nolimits}

\def \R{\mathop{\hbox{\main R}}\nolimits}
\def \S{\mathop{\hbox{\main S}}\nolimits}

\def \GL{\mathop{\hbox{\rm GL}}\nolimits}

\def \Rat{\mathop{\hbox{\rm Rat}}\nolimits}

\def \Res{\mathop{\hbox{\rm Res}}\nolimits}

\def \diag{\mathop{\hbox{\rm diag}}\nolimits}

\def \GL{\mathop{\hbox{\rm GL}}\nolimits}

\def \deg{\mathop{\hbox{\rm deg}}\nolimits}

\def \Stab{\mathop{\hbox{\rm Stab}}\nolimits}

\topmatter
\title Une remarque sur le degr\'e formel d'une s\'erie
discr\`ete d'un groupe lin\'eaire g\'en\'eral $p$-adique \endtitle

\rightheadtext{Une remarque sur le degr\'e formel d'une s\'erie
discr\`ete}
\author Volker Heiermann \endauthor

\address Institut f\"ur Mathematik, Humboldt-Universit\"at zu
Berlin, Unter den Linden 6, 10099 Berlin, Allemagne \endaddress

\email heierman\@mathematik.hu-berlin.de \endemail

\abstract We show by the example of the general linear group, how
one can get from \cite {2} precise information on the formal
degree of a square integrable representation of a $p$-adic group.

\null \null \noindent{R\'ESUM\'E: Nous montrons \`a l'exemple du
groupe lin\'eaire g\'en\'eral, comment on peut d\'eduire de
\cite{2} des informations pr\'ecises sur le degr\'e formel d'une
repr\'esentation de carr\'e int\'egrable d'un groupe $p$-adique.}
\endabstract

\endtopmatter
\document

\null \null Soit $F$ un corps local non archim\'edien de valeur
absolue normalis\'ee $\vert .\vert =\vert .\vert _F$ et $m>0$ un
entier. Fixons une repr\'esentation cuspidale unitaire $\sigma $
de $\GL _m(F)$ et un entier $d\geq 1$. Posons $n=md$, $G=\GL_n
(F)$. Notons $M$ l'unique sous-groupe de Levi standard de $G$ qui
est isomorphe \`a $\GL _m(F)\times\cdots\times\GL _m(F)$, $\det
_m$ le d\'eterminant de $\GL _m(F)$ et $\pi _d$ l'unique s\'erie
discr\`ete de $\GL _n(F)$ qui est un sous-quotient de l'induite
parabolique (normalis\'ee) de $\sigma\vert\det _m\vert^{(d-1)/2}$
$\otimes\cdots \otimes\sigma\vert\det _m\vert^{(-d+1)/2}:=\rho
_d$. On va appliquer les r\'esultats de \cite{2} pour calculer le
degr\'e formel de $\pi _d$ en fonction de celui de $\sigma $. On
retrouvera \`a cette occasion le r\'esultat de A.-M. Aubert et R.
Plymen \cite{1} dont la preuve utilisait la th\'eorie des types de
Bushnell-Kutzko.

\null {\bf 1. Les mesures}

\null Notons $\alpha _1,\dots ,\alpha _{d-1}$ les racines simples
de $G$ qui ne sont pas des racines de $M$ et identifions-les \`a
des racines du tore d\'eploy\'e maximal de $M$. Lorsque $l$ est un
entier, $1\leq l\leq d$, $M_l$ d\'esignera le sous-groupe de Levi
standard contenant $M$, obtenu en adjoignant \`a $M$ les racines
$\alpha _l,\dots ,\alpha _{d-1}$. Il est isomorphe \`a $\GL
_m(F)\times \cdots\times\GL _m(F)\times\GL _{(d-l+1)m}(F)$. En
particulier $M_1=G$ et $M_d=M$. Nous noterons $(m_1,\dots ,m_l)$
un \'el\'ement g\'en\'eral de $M_l$ (o\`u $m_1,\dots
,m_{l-1}\in\GL _m(F)$ et $m_l\in\GL _{(d-l+1)m}(F)$). Le symb\^ole
$M_l^1$ d\'esignera l'intersection des noyaux des caract\'eres non
ramifi\'es de $M_l$. Le g\'en\'erateur $>1$ de l'image de $\vert
.\vert$ sera not\'e $q$.

\null Notons $\X (M_l)$ le groupe des caract\`eres non ramifi\'es
de $M_l$ et $\Xim (M_l)$ le sous-groupe form\'e des caract\`eres
unitaires. Identifions le groupe $\Xim (M_l)$ \`a $(S^1)^l$ au
moyen de l'isomorphisme qui envoie un \'el\'ement $(q^{s_1},\cdots
,q^{s_l})$ de $(S^1)^l$ sur le caract\`ere non ramifi\'e
$\vert\det _m\vert_F^{s_1}$ $\cdots\vert\det _m\vert_F^{s_{l-1}}$
$\vert\det _{(d-l)m}\vert_F^{s_l}$ de $M_l$.

\null Le tore d\'eploy\'e maximal dans le centre de $M_l$ sera
not\'e $T_{M_l}$. Il est isomorphe \`a $(F^{\times })^l$. Le
groupe $\Xim (T_{M_l})$ est isomorphe \`a $(S^1)^l$ au moyen de
l'isomorphisme qui envoie un \'el\'ement $(q^{s_1},\cdots
,q^{s_l})$  de $(S^1)^l$ sur le caract\`ere non ramifi\'e de
$(F^{\times })^l$ donn\'e par $\vert\det _1\vert_F^{s_1}\cdots
\vert\det _1\vert_F^{s_l}$.

\null Suivant \cite{2} on munit le groupe $\Xim (T_{M_l})$ de
l'unique mesure de Haar de mesure totale \'egale \`a $1$. La
mesure sur $\Xim (M)$ est l'unique mesure de Haar telle que la
restriction $\Xim (M) \rightarrow \Xim (T_M)$ pr\'eserve
localement les mesures. En identifiant les deux groupes avec le
tore $(S^1)^l$, cette restriction correspond \`a l'application de
$(S^1)^l$ dans $(S^1)^l$ qui envoie $(z_1,\dots ,z_l)$ sur
$(z_1^m,\dots , z_{l-1}^m, z_l^{(d-l+1)m})$. Il en r\'esulte que
$\Xim (M)$ a la mesure $(d-l+1)m^l$.

\null Notons $\o _l$ l'orbite inertielle de $\rho :=\sigma \otimes
\sigma\otimes\dots\otimes\sigma$ par $\X (M_l)$, $\o _{l,0}$
l'orbite unitaire de $\rho $ par $\Xim (M_l)$ et $t$ l'ordre du
stabilisateur de $\sigma $ dans $\Xim (M)$. La mesure sur $\o
_{l,0}$ est l'unique mesure de Haar telle que l'application $\Xim
(M_l)\rightarrow\o _{l,0}$ qui envoie un caract\`ere unitaire non
ramifi\'e $\chi $ sur (la classe d'\'equivalence de)
$\sigma\otimes\chi $ pr\'eserve localement les mesures. Comme les
fibres de cette application ont tous m\^eme cardinalit\'e $t^l$,
la mesure de $\o _{l,0}$ est $(d-l+1)(m/t)^l$.

\null Tout sous-groupe ferm\'e $H$ de $G$ sera muni de l'unique
mesure de Haar telle que son intersection avec $K=\GL _n(O_F)$ a
la mesure $1$, $O_F$ d\'esignant l'anneau des entiers de $F$.

\null {\bf 2. La racine $\ti {\alpha }_{l-1}$}

\null Consid\'erons $\alpha _{l-1}$ comme racine de $T_{M_l}$.
Fixons un g\'en\'erateur $\ti{\omega }$ de l'id\'eal maximal de
$F$. Un g\'en\'erateur $h_{l-1}$ du $\Bbb Z$-module $M_{l-1}^1\cap
M_l/M_l^1$ est donn\'e par la classe modulo $M_l^1$ d'une matrice
diagonale $\diag (x_1,\dots ,x_n)$ avec $x_{m(l-1)}=\ti{\omega }$,
$x_{m(l-1)+1}=\ti{\omega }^{-1}$ et $x_i=1$ sinon. Notons $\Rat
(M_l)$ le groupe des caract\`eres alg\'ebriques de $M_l$ d\'efinis
sur $F$.

\null L'\'el\'ement de $\Rat (M_l)$ qui correspond \`a
$(d-l+1)m\alpha _{l-1}$ envoie $(m_1,\dots , $ $m_l)$ sur  $\det
_m(m_{l-1})^{d-l+1}\det _{(d-l+1)m}(m_l)^{-1}$. L'\'el\'ement
$\alpha _{l-1} ^*:=(d-l+1/d-l+2)m\alpha _{l-1}$ v\'erifie donc
$\langle\alpha _{l-1}^*, H_M$ $(h_{l-1})\rangle $ $=1$. La
puissance $h_{l-1}^t$ de $h_{l-1}$ est minimale telle que $\chi\in
\X(M_l)$ v\'erifie $\chi (h_{l-1}^t)=1$, si et seulement si
$\chi\in \X(M_{l-1})\Stab (\rho )$. On trouve donc $\ti {\alpha
}_{l-1}=(d-l+1/d-l+2)(m/t)\alpha _{l-1}$ dans les notations de
[2].

\null Notons $\alpha _l^{\vee }$ la coracine associ\'ee \`a
$\alpha _l$. Un calcul \'el\'ementaire donne $t\ \langle\alpha
_l^{\vee },z_1\ti{\alpha }_1+z_2\ti{\alpha }_2+\cdots
+z_{d-1}\ti{\alpha }_{d-1}\rangle =z_l-{d-l-1\over d-l}z_{l+1}$
(o\`u $z_d:=0$). Par suite, avec $r_l=t {d-l+1\over 2}$, on trouve
$\rho _d=\rho\otimes\chi _{r_1\ti{\alpha }_1+r_2\ti{\alpha
}_2+\cdots +r_{d-1}\ti {\alpha }_{d-1}}$.

\null {\bf 3. La fonction $\mu $}

\null Notons $a=n(\sigma\times\sigma ^{\vee })$ le conducteur
d'Artin de paires (cf. \cite{3}, p.291) et $\mu ^{M_l}$ la
fonction $\mu $ de Harish-Chandra d\'efinie sur $\o:=\o _d$
relative \`a $M_l$ comme produit des op\'erateurs d'entrelacement
(cf. \cite {2} {\bf 1.5}), $\mu :=\mu ^G$. Posons
$a_{M_l}^*=\Rat(M_l)\otimes _{\Bbb Z}\Bbb R$. On a une surjection
canonique $a_{M_l,\Bbb C}^*\rightarrow\X (M_l)$,
$\lambda\mapsto\chi _{\lambda }$. Le r\'esultat suivant est une
cons\'equence imm\'ediate de la formule de produit pour $\mu
^{M_l}$, utilisant la formule explicite pour $\mu $ dans le cas
d'un sous-groupe de Levi maximal (cf. \cite {1}, thm. {\bf 3.3}).

\null {\bf 3.1 Proposition:} \it Posons $\alpha _{i,j}^{\vee
}=\alpha _i^{\vee }+\cdots +\alpha _j^{\vee }$ pour $i\leq j$. La
fonction $\lambda\mapsto\mu(\rho\otimes\chi _{\lambda })$ est
r\'eguli\`ere et non nulle en dehors de la r\'eunion des
hyperplans affines de $a_M^*$ de la forme $\langle\alpha
_{i,j}^{\vee },\lambda\rangle=0,\pm 1$ avec $1\leq i<j\leq d-1$.

La fonction $({\mu^{M_{l-1}}/ \mu ^{M_l}})(\rho _0\otimes\chi
_{z\ti{\alpha }_{l-1}+r_l\ti{\alpha }_l+\dots +r_{d-1}\ti {\alpha
}_{d-1}})$ vaut
$$q^{(d-l+1)a}{(1-q^{t{d-l\over 2}-z})(1-q^{t{d-l\over 2}+z})\over
(1-q^{-t{d-l\over 2}-t-z})(1-q^{-t{d-l\over 2}-t+z})}.$$ \rm

\null  {\bf 4. La donn\'ee de r\'esidu $\Res _A$}

\null Posons $A_l=\rho_d\X (M_l)$, notons $r(A_l)$ "l'origine" de
$A_l$ (cf. \cite{2} {\bf 1.4}), $A_{l,0}$ le sous-espace de $A_l$,
form\'e des points de partie r\'eelle $r(A_l)$, et $\S _{A_1}$
l'ensemble form\'e des hyperplans affines $\{\rho\otimes\chi
_{\lambda }\vert \langle\alpha _l^{\vee },\lambda\rangle =1\}$,
$l=1,\dots, d-1$, de $\o :=\o _d$. D\'esignons par $\R(\S _{A_1})$
l'espace des fonctions rationnelles sur $\o _d$, r\'eguli\`eres en
dehors des hyperplans affines dans $\S _{A_1}$. Remarquons que la
partie r\'eelle $\Re (\rho\otimes\chi _{\lambda }):=\Re (\lambda
)$ est bien d\'efinie. Le symb\^ole $\int _{\Re (\rho ')=R}
d_{A_d}\Im (\rho ')$ d\'esignera la mesure sur $\chi_R\o _0$,
d\'eduite de celle sur $\o _0=A_{d,0}$. De fa\c con analogue pour
$\int _{A_{l,0}} d_{A_l}\Im(\rho ')$. L'ordre sur $a_M^*$ induit
par $P$ sera not\'e $>_P$.

\null {\bf 4.1 Proposition:} \it Soit $\psi $ dans $\R(\S _{A_1})$
invariante par $\X (G)$. Pour $\chi\in\X (G)$ et $z_1,\dots
,z_{d-1}\in \Bbb C$, posons $f(z_1,\cdots ,z_{d-1})=\psi
(\rho\otimes\chi\chi _{z_1\ti{\alpha }_1+\dots +z_{d-1}\ti {\alpha
}_{d-1}})$. On a
$$\int _{\Re (\rho ')=R\gg _P0}\psi (\rho ') d_{A_d}\Im(\rho ')=\sum _{l=1}^d\int
_{A_{l,0}}(\Res _{A_l}\psi )(\rho ') d_{A_l}\Im(\rho '),$$ avec
$\Res _{A_l}\psi (\rho\otimes\chi\chi _{z_1\ti{\alpha }_1+\dots +
z_{l-1}\ti {\alpha }_{l-1}})$ \'egal \`a
$$({m\log q\over t})^{d-l}{1\over d-l+1}\Res_{z_l=r_l}
(\dots(\Res_{z_{d-1}=r_{d-1}}f).)(z_1,\dots ,z_{l-1}).$$

\null Preuve: \rm \'Ecrivons $R=R_1\ti{\alpha }_1+R_2\ti{\alpha }
_2+\cdots +R_{d-1}\ti {\alpha }_{d-1}$. Par la suite exacte dans
\cite{2} {\bf 1.2} et notre choix des mesures, on a
$$\eqalign {&{\log q\over 2\pi } \int _{\Re (\rho ')=R}\psi (\rho ')d_{A_d}\Im(\rho ') \cr
=&({\log q\over 2\pi })^{d}({m\over t})^d\int _0^{2\pi\over\log
q}\dots\int _0^{2\pi\over\log q} f(R_1+it_1,\dots
,R_{d-1}+it_{d-1})dt_{d-1}\dots dt_1.\cr}$$ Si on fixe $z_1,\dots
,z_{l-1}$ avec $\Re (z_i)=R_i$, la fonction $z_l\mapsto\psi
(z_1\ti{\alpha }_1+\cdots +z_{l-1}\ti{\alpha }_{l-1}+z_l\ti{\alpha
_l}+r_{l+1}z_{l+1}+\cdots r_{d-1}z_{d-1})$ n'a, compte tenue de la
proposition {\bf 3.1} et du calcul dans {\bf 2.}, au plus un
p\^ole en $z_l=r_l$. Par les arguments dans [2] {\bf 3.6},
l'int\'egrale vaut
$$\eqalign {&({1\over 2\pi })^l ({m\log q\over t})^d\sum _{l=1}^d\int
_0^{2\pi\over\log q}\dots\int _0^{2\pi\over\log q} \cr &\qquad\Res
_{z_l=r_l}(\dots (\Res _{z_{d-1}=r_{d-1}} f).) (it_1,\dots
,it_{l-1})dt_{l-1}\dots dt_1\cr =&\sum _{l=1}^d({m\log q\over
t})^{d-l}{1\over d-l+1} (d-l+1) ({1\over 2\pi })^l ({m\log q\over
t})^l \cr &\qquad \int _0^{2\pi\over\log q}\dots \int
_0^{2\pi\over\log q}\Res _{z_l=r_l}(\dots (\Res
_{z_{d-1}=r_{d-1}}f).) (it_1,\dots ,it_{l-1})dt_{l-1}\dots dt_1
\cr =&{\log q\over 2\pi}\sum _{l=1}^d \int _{A_{l,0}} (\Res
_{A_l}\psi )(\rho ')d_{A_l}\Im (\rho '). }$$

\null {\bf 4.2 Corollaire:} \it Avec $f(z_1,\cdots ,z_{d-1})=\mu
(\rho\otimes\chi _{z_1\ti{\alpha }_1+\dots +z_{d-1}\ti {\alpha
}_{d-1}})$, on a
$$(\Res _{A_1}^P\mu )=({m\log q\over t})^{d-1}{1\over
d}\Res_{z_1=r_1}(\dots(\Res_{z_{d-1}=r_{d-1}}f).).$$

\null Preuve: \rm Par [2] {\bf 3.9}, la donn\'ee de r\'esidu $\Res
_{A_1}^P$ est d\'etermin\'e par sa restriction \`a $\S _{A_1}$,
donc \'egale \`a l'expression donn\'ee dans la proposition {\bf
4.1}. \hfill {\fin 2}

\null {\bf 5. Le degr\'e formel}

\null {\bf 5.1 Th\'eor\`eme:} \it Le degr\'e formel de l'unique
sous-quotient irr\'eductible $\pi _d$ de la repr\'esentation
induite parabolique (normalis\'ee) de la repr\'esentation
$\sigma\vert\det _m\vert^{(d-1)/2}$ $\otimes\cdots
\otimes\sigma\vert\det _m\vert^{(-d+1)/2}$ de $M$ est li\'e \`a
celui de $\sigma $ par la formule
$$\deg (\pi _d)={\vert\GL _n(\hbox{\ens F}_q)\vert\over\vert\GL _m(\hbox{\ens F}_q)\vert^d}
\ q^{mn-n^2}\ {m^{d-1}\over t^{d-1}d}\ q^{a{d(d-1)\over 2}}\ q^{td
(d-1)\over 2}\ {(q^t-1)^d\over q^{td}-1}\ \deg(\sigma )^d.$$

\null Preuve: \rm Par la remarque dans le paragraphe {\bf 8.6} de
\cite{2}, on a
$$\deg (\pi _d)=\gamma (G/M) \deg (\rho _d) \vert\Stab
(A_1)\vert^{-1} (\Res _{A_1}^P\mu )(\rho _d).$$ Comme $\rho _d$
est r\'egulier, $\Stab (A_1)=\{1\}$. Il est imm\'ediate par
d\'efinition du degr\'e formel que celui d'un produit tensoriel de
repr\'esentations est \'egal au produit des degr\'es formels. La
constante $\gamma (G/M)$ est \'egal \`a ${\vert\GL _n(\hbox{\ens
F}_q)\vert\over\vert\GL _m(\hbox{\ens F}_q)\vert^d}q^{mn-n^2}$,
comme on le d\'eduit directement de la formule pour $\gamma (G/M)$
dans \cite{4} en haut de la page 241, en posant $H=I_n+\ti{\omega
}\GL _n(O_F)$. Il reste \`a calculer le r\'esidu, utilisant {\bf
3.1}.
$$\eqalign {&(\Res _{A_1}^P\mu )(\rho _d)\cr
=&({m\log q\over t})^{d-1} {1\over d} \cr &\qquad \prod _{l=2}^d\
q^{a(d-l+1)} {(1-q^{t{d-l\over 2}-r_{l-1}}) (1-q^{t{d-l\over
2}+r_{l-1}})\over (1-q^{-t{d-l\over 2}-t-r_{l-1}})}\ \Res
_{z=r_{l-1}} {1\over 1-q^{-t{d-l\over 2}-t+z}}\cr =&({m\log q\over
t})^{d-1} {1\over d}\ q^{a{d(d-1)\over 2}}\prod _{l=2}^d\
{(q^{-t}-1) (q^{t(d-l+1)}-1)\over (q^{-t(d-l+2)}-1)}\ \Res
_{z=r_{l-1}} {1\over q^{-t{d-l\over 2}-t+z}-1}\cr =&({m\log q\over
t})^{d-1} {1\over d}\ q^{a{d(d-1)\over 2}}(q^{-t}-1)^{d-1}\
{q^t-1\over q^{-dt}-1}({1\over \log q})^{d-1} \prod _{l=2}^{d-1}
{q^{t(d-l+1)}-1\over q^{-t(d-l+1)}-1}\cr =&({m\over t})^{d-1}
{1\over d}\ q^{a{d(d-1)\over 2}}(q^{-t}-1)^{d-1}{q^t-1\over
q^{-dt}-1} \prod _{l=2}^{d-1} (-q^{t(d-l+1)})\cr =&({m\over
t})^{d-1} {1\over d}\ q^{a{d(d-1)\over 2}} q^{-t(d-1)}(-1)^{d-1}
(q^t-1)^{d-1} (-q^{td}) {q^t-1\over q^{td}-1} (-1)^{d-2}\prod
_{l=2}^{d-1}q^{tl}\cr =&({m\over t})^{d-1} {1\over d}\
q^{a{d(d-1)\over 2}} {(q^t-1)^{d}\over q^{td}-1}\ q^{t{d(d-1)\over
2}}.\cr}$$ \srm L'auteur remercie l'Universit\'e Purdue ainsi que
F. Shahidi pour leur hospitalit\'e. \rm


\Refs
\ref \key{1} \by A.-M. Aubert, R. Plymen \paper Plancherel
measure for {\rm GL}$(n)$: Explicit Formulas and Bernstein
Decomposition \jour preprint, \yr 2003\endref

\ref \key{2} \by V. Heiermann \paper D\'ecomposition spectrale
d'un groupe r\'eductif $p$-adique \jour \`a para\^\i tre dans
Journal de l'Institut de Math\'ematiques de Jussieu
\endref

\ref \key{3} \by F. Shahidi \paper Langlands' conjecture on
Plancherel measures for $p$-adic groups \jour Progr. Math. \vol
101 \publ Birkh\"auser \publaddr Boston, MA \yr 1991 \pages
277--295\endref

\ref \key{4} \by J.-L. Waldspurger \paper La formule de Plancherel
pour les groupes $p$-adiques (d'apr\`es Harish-Chandra) \jour
Journal de l'Institut de Math\'ematiques de Jussieu, \vol 2 \yr
2003 \pages 235--333\endref
\endRefs
\enddocument